\documentclass[pdflatex]{sn-jnl}
\usepackage[T1]{fontenc}
\usepackage{lmodern}
\usepackage{newlfont}
\usepackage{dsfont}
\usepackage{amsmath}
\usepackage{amsfonts}
\usepackage{amssymb}

\usepackage{xcolor}
\usepackage{interval}
\intervalconfig{soft open fences}
\usepackage{enumitem}
\setlist[enumerate]{labelindent=\parindent,leftmargin=*}

\usepackage[numbers,sort&compress]{natbib}
\usepackage{graphicx}
\newtheorem{thm}{Theorem}[section]
\newtheorem{cor}[thm]{Corollary}

\newtheorem{prop}[thm]{Proposition}
\theoremstyle{definition}
\newtheorem{defn}[thm]{Definition}
\theoremstyle{remark}
\newtheorem{rem}[thm]{Remark}
\numberwithin{equation}{section}

\newcommand{\abs}[1]{\left\vert#1\right\vert}
\newcommand{\set}[1]{\left\{#1\right\}}

\newcommand{\C}{\mathbb{C}}
\newcommand{\Cn}{\mathbb{C}^n}
\newcommand{\D}{\mathbb{D}}

\begin{document}
    \title{A Canonical Characterization of Normal Functions}

    \author[1]{\fnm{Peter V.} \sur{Dovbush}}
    \email{peter.dovbush@gmail.com}
    \author[2]{\fnm{Steven G.} \sur{Krantz}}
    \email{sk@math.wustl.edu}

    \affil[1]{\orgdiv{Institute of Mathematics and Computer Science}, \orgname{Moldova State University}, \orgaddress{\street{5 Academy Street}, \city{Kishinev}, \postcode{MD-2028}, \country{Republic of Moldova}}}

    \affil[2]{\orgname{Washington University in St.\ Louis}, \orgaddress{\street{Street}, \city{St.\ Louis}, \postcode{63130}, \state{Missouri}, \country{USA}}}


    \abstract{ We characterize normal families in the unit ball as those families of analytic functions whose restrictions to each complex line through the origin are normal.
    We then generalize this result to a characterization of normal functions according to behavior on analytic disks.
    We give a simple proof of an old theorem of Hartogs stating that a formal power series at $0$ in $\Cn$ is convergent if its restriction to each complex line through the origin is convergent.}


    \maketitle

    \section{Introduction}

    Let $\D:=\set{z\in \C \mid \abs{z}<1}$ be the unit disk in $\C$.
    Then the Poincar\'{e} metric
    \[
        d \rho^2=\frac{2 d z d \bar{z}}{\left(1-|z|^2\right)^2}\,.
    \]

    Let $\Omega$ be a domain (connected open set) in $\C^n$.
    Denote by $\mathcal{O}(\Omega)$ the ring of all holomorphic functions on $\Omega$.Let $\D(\Omega)$ be the set of holomorphic mappings from $\Omega$ to $\D$, and $\Omega(\D)$ the set of holomorphic mappings from $\D$ to $\Omega$.

    A family $\mathcal{F}\subset \mathcal{O}(\Omega)$ of holomorphic functions on a domain $\Omega \subset \Cn$ is said to be \emph{normal in $\Omega$ in the sense of Montel} if every sequence $\{f_j\} \subset \mathcal{F}$ contains either a subsequence which converges uniformly on every compact subset of $\Omega$ to a finite limit function, or a subsequence which converges uniformly to $\infty$ on every compact subset of $\Omega$.

    A family $\mathcal{F}\subset \mathcal{O}(\Omega)$ is said to be normal at a point $z_0 \in \Omega$ if it is normal in some neighborhood of $z_0$.

    Suppose that $f(z)$ is meromorphic in the unit disk $\D$.
    Consider the family $\{f_a(z)\}$ consisting of all functions of the form
    \[
        f_a(z) \equiv f\!\left(\frac{z - a}{\bar{a}z - 1}\right)\,,
    \]
    where $a$ ranges over the unit disk $\D$.

    Following K.~Yosida~\cite{zbMATH03013560}, we say that the function $f(z)$ belongs to the class $(A)$ if the family $\{f_a(z)\}$ is a normal family in the unit disk $\D$ in the sense of Montel.

    We have

    \begin{thm}
        \label{thm:yosida-noshiro}
        A function $f$ belongs to the class $(A)$ if and only if exist a positive number $C$ such that
        \[
            \frac{\left(1-|z|^2\right)|f'(z)|}{1+|f(z)|^2}<C
        \]
        for $|z|<1$.
    \end{thm}

    \begin{rem}
        The above theorem corresponds to Yosida's~\cite[Theorem 1]{zbMATH03013560}, and was proved in Noshiro~\cite[Theorem 1]{zbMATH03034393}.
        It is worth noting that we are invoking here the spherical metric on the range.
    \end{rem}

    Let $U$ and $V$ be open subsets of $\Cn$.
    A holomorphic map $f: U \rightarrow V$ that is one-to-one and onto is called a \emph{biholomorphism} from $U$ to $V$.
    A biholomorphism from $U$ to $U$ is called an \emph{automorphism} of $U$.
    The set of all automorphisms of $U$ is denoted by $\operatorname{Aut}(U)$ is a group under composition.

    For every $a \in \mathbb{D}$,
    \[
        \varphi_a(z)=\frac{a-z}{1-\bar{a} z}
    \]
    is an automorphism of $\mathbb{D}$.

    One can show by Schwarz's lemma that a map $\varphi_a$ is an automorphism of $\mathbb{D}$ if and only if there are $\theta \in \mathbb{R}$ and $a \in \mathbb{D}$ such that
    \[
        \varphi(z)=e^{i \theta}\cdot \varphi_a(z)\,.
    \]

    Let $f(z)$ be meromorphic in $\D$.
    For any real $\theta$ and $a$, such that $|a|<1$, we shall call
    \[
        f(z, a, \theta)=f\left(e^{i \theta}\cdot \frac{a-z}{1-\bar{a} z}\right)
    \]
    a \emph{translate} of $f(z)$.

    A family $\mathcal{F}$ of functions $f(z)$ is said to be invariant (Hayman~\cite{MR0070716}) if whenever $f(z) \in \mathcal{F}$ then all the translates $f(z, a, \theta) \in \mathcal{F}$.
    Many normal families turn out also to be invariant, and for these it is possible to obtain bounds of considerably greater precision than for general normal families.

    Following Lehto and Virtanen~\cite{MR0087746} a function $f(z)$ is called normal if $f(z)$ together with its translates constitutes a normal family.
    Normal meromorphic functions admit the following characterization in terms of the spherical derivative:

    \begin{thm}
        \label{thm:criterion-norm}
        An invariant family $\mathcal{F}$ of functions $f(z)$ meromorphic in $\D$ is normal there if and only if there is a constant $B$, such that
        \[
            \frac{|f'(0)|}{1+|f(0)|^2}\leqslant B
        \]
        for $f(z) \in F$.
        A function $f(z)$ meromorphic in $\D$ is normal if and only if
        \begin{equation}
            \frac{\left(1-|z|^2\right)|f'(z)|}{1+|f(z)|^2}\leqslant B \quad(|z|<1)\,,
            \label{eq:basic}
        \end{equation}
        where $B$ is a constant.
    \end{thm}

    The inequality~\eqref{eq:basic} can be rewritten in the form
    \begin{equation}
        \frac{|f'(z)|^2 dzd\overline{z}}{\left(1+|f(z)|^2\right)^2}\leqslant B\frac{2dzd\overline{z}}{\left(1-|z|^2\right)^2}\quad(|z|<1)\,,
        \label{eq:basic2}
    \end{equation}

    Remark that in the right side of inequality~\eqref{eq:basic2} stands the Poincar\`{e} metric tensor.

    \section{Invariant Families in the Unit Ball}

    For the unit ball $\mathbb{B}=\{z\in \Cn \mid |z|<1\}$ the Bergman metric is defined by the formula
    \[
        d s^2=(n+1)\left\{\frac{|d z|^2}{1-|z|^2}+\sum_{\mu, \nu=1}^n \frac{\bar{z}_\mu z_\nu d z_\mu d \bar{z}_\nu}{\left(1-|z|^2\right)^2}\right\}\,,
    \]
    where $|d z|^2=\sum_{\nu=1}^n|d z_\nu|^2$.

    For $n=1$ and $r=1$, the Bergman metric coincides with the Poincar\'e metric in the unit disk $\D$.

    The interested reader can refer to~\cite{MR3114665} for further literature on the Bergman kernel and metric.

    For every function $\varphi$ of class $C^2(\Omega)$ we define at each point $z\in \Omega$ a hermitian form
    \[
        L_z(\varphi, v):=\sum_{k,l=1}^n \frac{\partial^2\varphi}{\partial z_k \partial \overline{z}_l}(z) v_k \overline{v}_l\,,
    \]
    and call it the \emph{Levi form} of the function $\varphi$ at $z$.
    Since $\varphi$ is real this form is Hermitian $\left(\partial^2 \varphi / \partial z_\mu \partial \bar{z}_\nu=\bar{\partial}^2 \varphi / \partial z_\nu \partial \bar{z}_\mu\right)$ and hence it takes real values on all vectors $v \in \Cn$.

    If $\psi:U\rightarrow\C^m$ is a holomorphic mapping and $\varphi\in C^2$ in a neighborhood of $\psi(z)$, then
    \[
        L_z(\varphi\circ\psi,v) = L_{\psi(z)}(\varphi,\psi_{*}v)\,,
    \]
    where $\psi_{*}=d \psi$ is the differential of the mapping $\psi$ at $z$.
    The last property shows that the Levi form is invariant relative to biholomorphic mappings.

    For a holomorphic function $f$ in $\Omega$, set
    \begin{equation}
        f^\sharp (z):=\sup_{ |v|=1}\sqrt{L_z(\log(1+|f|^2), v)}\,.
        \label{eq:e2}
    \end{equation}
    This quantity is well defined since the Levi form $L_z(\log(1+|f|^2), v)$ is nonnegative for all $z\in \Omega$.

    For all $z\in \Omega$ and $v\in \mathbb{C}^n$
    \[
        L_z\left(\log \left(1+|f|^2\right), v\right)\equiv \frac{|d f_z(v)|^2}{\left(1+|f(z)|^2\right)^2}\,.
    \]
    In particular, for $n = 1$ the formula~\eqref{eq:e2} takes the form
    \begin{equation}
        f^\sharp (z)=\frac{|f'(z)|}{1+|f(z)|^2}\,.
        \label{eq:e}
    \end{equation}
    The quantity $f^\sharp (z)$ is called \emph{spherical derivative } of $f$ at a the point $z$.

    The following generalization of Marty's theorem was proved in~\cite[Theorem 2.1]{MR4071476}.

    \begin{thm}
        \label{thm:marty:1}
        A family $\mathcal{F}$ of functions holomorphic on $\Omega$ is normal on $\Omega\subset \Cn$ if and only if for each compact subset $K\subset \Omega$ there exists a constant $M(K)$ such that at each point $z\in K$
        \begin{equation}
            f^\sharp(z)\leq M(K)
            \label{eq:marty:1}
        \end{equation}
        for all $f\in \mathcal{F}$.
    \end{thm}

    \begin{defn}
        \label{defn:norm-manif-1}
        We say that a holomorphic function $f: \mathbb{B}\rightarrow \C$ is \emph{normal} if the family
        \[
            {\mathcal{F}}= \set{f \circ g \mid g \in \operatorname{Aut}(\mathbb{B})}
        \]
        is normal.
    \end{defn}

    Repeating the proof of Theorem~\ref{thm:criterion-norm} given in Noshiro~\cite[p.~87]{MR0133464}, one obtains at once the following:

    \begin{prop}
        \label{prop:norm-manif-1}
        Let $ds^2$ be the Bergman metric on $\mathbb{B}$ in $\C^n$.
        If a holomorphic function $f: \mathbb{B}\rightarrow \C$ satisfies
        \[
            L_z(\log(1+|f|^2), v)\leq C\cdot ds^2(z,v)
        \]
        for a finite constant $C$, then $f$ is a normal holomorphic function.
    \end{prop}

    \begin{proof}
        Let
        \[
            \mathcal{F}=\{f\circ g:g\in\operatorname{Aut}(\mathbb{B})\}\,.
        \]
        Since the Bergman metric is $\operatorname{Aut}(\mathbb{B})$-invariant, for every $g\in\operatorname{Aut}(\mathbb{B})$ we have
        \[
            \begin{aligned}
                L_z\left(\log(1+|f\circ g|^2),v\right) & = L_{g(z)}\left(\log(1+|f|^2),g_{*}v\right) \\
                                                           & \leq C\,ds^2(g(z),g_{*}v)                   \\
                                                           & = C\,ds^2(z,v)\,.
            \end{aligned}
        \]
        Consequently, on every compact subset $K\Subset\mathbb{B}$, the spherical derivatives of the functions in $\mathcal{F}$ are uniformly bounded.
        Indeed, on $K$ the Bergman metric is uniformly comparable with the Euclidean metric.
        Thus $\mathcal{F}$ is locally equicontinuous with respect to the spherical distance on $\overline{\C}$.

        Since $\overline{\C}$ is compact, the Arzel\`{a}--Ascoli theorem implies that $\mathcal{F}$ is a normal family.
        Hence $f$ is normal.
    \end{proof}

    \begin{defn}
        \label{defn:norm-manif-2}
        We say that a subset ${\mathcal{F}}$ of ${\mathcal{O}}(\mathbb{B}, \C)$ is \emph{$\operatorname{Aut}(\mathbb{B})$-invariant} if $f \circ g \in{\mathcal{F}}$ for every $f \in{\mathcal{F}}$ and every $g \in \operatorname{Aut}(\mathbb{B})$.
    \end{defn}

    \begin{thm}
        \label{thm:norm-manif-1}
        If a subset $\mathcal{F}$ of $\mathcal{O}(\mathbb{B}, \C)$ is an $\operatorname{Aut}(\mathbb{B})$-invariant normal family, then there exists a constant $C$ such that
        \begin{equation}
            L_z(\log(1+|f|^2),v)\leq C \cdot ds^2(z,v)
            \label{eq:norm-inequality}
        \end{equation}
        for every $f \in \mathcal{F}$.
    \end{thm}

    \begin{proof}
        For $z\in\mathbb{B}$, define
        \[
            C(z) := \sup_{\substack{f\in\mathcal{F}\\v\ne0}}\frac{ L_z\left(\log(1+|f|^2),v\right) }{ ds^2(z,v) }\,.
        \]

        We first show that $C(0)<\infty$.

        Since $\mathcal{F}$ is a normal family, Marty's theorem (Theorem~\ref{thm:marty:1}) implies that the spherical derivatives of the functions in $\mathcal{F}$ are locally uniformly bounded.
        In particular, there exists a constant $M<\infty$ such that
        \[
            L_0\left(\log(1+|f|^2),v\right)\leq M|v|^2
        \]
        for every $f\in\mathcal{F}$ and every $v\in\C^n$.

        We have
        \[
            2|v|^2=ds^2(0,v)
        \]
        for every $v\in\C^n$.
        Consequently,
        \[
            L_0\left(\log(1+|f|^2),v\right) \leq M|v|^2\leq Mds^2(0,v)\,.
        \]
        Hence
        \[
            C(0)\leq M<\infty\,.
        \]

        We next show that $C$ is invariant under $\operatorname{Aut}(\mathbb{B})$.
        Let $g\in\operatorname{Aut}(\mathbb{B})$.
        Since $\mathcal{F}$ is $\operatorname{Aut}(\mathbb{B})$-invariant, $f\circ g\in\mathcal{F}$ whenever $f\in\mathcal{F}$.
        By the chain rule,
        \[
            L_z\left(\log(1+|f\circ g|^2),v\right) = L_{g(z)}\left(\log(1+|f|^2),g_{*}v\right)\,.
        \]
        By the definition of $C(g(z))$,
        \[
            L_{g(z)}\left(\log(1+|f|^2),g_{*}v\right) \leq C(g(z))\,ds^2(g(z),g_{*}v)\,.
        \]
        The $\operatorname{Aut}(\mathbb{B})$-invariance of the Bergman metric gives
        \[
            ds^2(g(z),g_{*}v)=ds^2(z,v)\,.
        \]
        Therefore
        \[
            L_z\left(\log(1+|f\circ g|^2),v\right) \leq C(g(z))\,ds^2(z,v)\,.
        \]
        Taking the supremum over $f\in\mathcal{F}$ and $v\neq0$ gives
        \[
            C(z)\leq C(g(z))\,.
        \]
        Replacing $g$ by $g^{-1}$ gives the reverse inequality.
        Hence
        \[
            C(z)=C(g(z))\,.
        \]
        Since $\mathbb{B}$ is homogeneous under $\operatorname{Aut}(\mathbb{B})$, $C(z)$ is constant on $\mathbb{B}$.
        Since $C(0)<\infty$, this constant is finite, and therefore
        \[
            L_z\left(\log(1+|f|^2),v\right) \leq C\,ds^2(z,v)
        \]
        for every $f\in\mathcal{F}$, every $z\in\mathbb{B}$, and every $v\in\C^n$.
    \end{proof}

    \section{Alexander's Theorem and Hartogs's Theorem}

    Alexander~\cite{MR0357850} proved the following characterization of normal families of holomorphic functions on the unit ball.
    It is a ``radial'' analog of a theorem of Nishino~\cite{MR0179384}, who proved that a family of holomorphic functions on a domain in $\C^2$ is normal if its restrictions to the coordinate lines are normal.

    Nishino proves a separate normality theorem not only for families of holomorphic functions but also for families of meromorphic functions.

    We acknowledge that there have been prior follow-ups to Nishino's paper, notably a direct generalization by Toshiaki Terada~\cite{MR0304698}; and an extension to mappings by Theodore J.~Barth~\cite{MR0374462}.
    Moreover, Nishino's results were anticipated by Renato Caccioppoli~\cite{zbMATH03012381}.

    We shall use the following form of Alexander's theorem.

    \begin{thm}
        \label{thm:alexander}
        Let $\mathcal{F}$ be a family of holomorphic functions on the unit ball $\mathbb{B}\subset\Cn$.
        Then $\mathcal{F}$ is normal on $\mathbb{B}$ if and only if, for every $c\in\C^n$ with $\lVert c\rVert=1$, the family
        \[
            \mathcal{F}_c := \set{\,\lambda\mapsto f(\lambda c):f\in\mathcal{F}\,}
        \]
        is normal on $\D$.

        Moreover, the following local version holds.
        If, for every $c\in\C^n$ with $\lVert c\rVert=1$, the family $\mathcal{F}_c$ is normal in some neighborhood of $0$ in $\D$ (the neighborhood may depend on $c$), then $\mathcal{F}$ is normal in some neighborhood of $0$ in $\mathbb{B}$.
    \end{thm}

    \begin{rem}
        \label{rem:alexander}
        The local assertion in Theorem~\ref{thm:alexander} is important when applying the theorem to formal power series.
        In that situation, the restriction of a formal power series to a fixed complex line is initially known to converge only in a neighborhood of the origin, and the radius of convergence may depend on the direction of the line.
    \end{rem}

    The preceding theorem immediately gives the following characterization of normal functions.

    \begin{thm}
        \label{thm:alexander-plus}
        Let $\mathcal{F}$ be a family of holomorphic functions on $\mathbb{B}$.
        \begin{enumerate}
            \item[(a)] A holomorphic function $f$ on $\mathbb{B}$ is normal if and only if its restriction to every complex line through the origin is normal.

            \item[(b)] The family $\mathcal{F}$ is normal on $\mathbb{B}$ if and only if its restriction to every complex line through the origin is normal.
        \end{enumerate}
    \end{thm}

    \begin{proof}
        For part~(a), consider the one-element family
        \[
            \mathcal{F}=\set{f}\,.
        \]
        By Theorem~\ref{thm:alexander}, $\mathcal{F}$ is normal on $\mathbb{B}$ if and only if, for every $c\in\C^n$ with $\lVert c\rVert=1$, the family
        \[
            \set{\lambda\mapsto f(\lambda c)}
        \]
        is normal on $\D$.
        Since this is a one-element family, this is equivalent to the normality of the function
        \[
            \lambda\mapsto f(\lambda c)
        \]
        on $\D$.
        Hence $f$ is normal on $\mathbb{B}$ if and only if its restriction to every complex line through the origin is normal.

        Part~(b) is exactly the first assertion of Theorem~\ref{thm:alexander}.
    \end{proof}

    \begin{rem}
        \label{rem:alexander-proof}
        The difficult implication in Theorem~\ref{thm:alexander} is the converse implication: normality of all the radial restrictions must imply normality of the original family.
        In particular, one cannot obtain this implication merely by applying the one-variable normality criterion separately for each direction $c$; the resulting constants must be shown to be uniform in $c$.
        This uniformity is the substance of Alexander's theorem.
    \end{rem}

    We now apply the local form of Alexander's theorem to obtain a proof of Hartogs's theorem concerning formal power series.

    Let
    \[
        \C\{z_1,\ldots,z_n\}
    \]
    denote the ring of convergent power series in $z=(z_1,\ldots,z_n)$ over $\C$, and let
    \[
        \C[[z_1,\ldots,z_n]]
    \]
    denote the ring of formal power series in $z_1,\ldots,z_n$ over $\C$.

    K.~Weierstra\ss{} based the theory of holomorphic functions on convergent power series.
    It is therefore important to have criteria which imply the convergence of a formal power series.
    The following classical theorem of Hartogs~\cite{MR1511365} gives such a criterion.

    \begin{cor}
        \label{cor:hartogs}
        Let
        \[
            F(z)=\sum_{\alpha\in\mathbb{N}^n}c_\alpha z^\alpha \in\C[[z_1,\ldots,z_n]]
        \]
        be a formal power series.
        If the restriction of $F$ to every complex line through the origin is convergent in a neighborhood of the origin, then $F$ is convergent in a neighborhood of the origin.
    \end{cor}

    \begin{proof}
        Put
        \[
            \mathcal{F}=\{f_m\}, \qquad f_m(z)=\sum_{|\alpha|\leq m}c_\alpha z^\alpha\,.
        \]
        By assumption, for every $c\in\C^n$ with $\lVert c\rVert=1$, the formal power series
        \[
            \sum_{\alpha\in\mathbb{N}^n}c_\alpha(\lambda c)^\alpha
        \]
        converges in some neighborhood of $0$ in $\D$.
        Hence there is a neighborhood $U$ of $0$ and a subsequence $\{f_{m_j}\}$ which converges locally uniformly on $U$ to a holomorphic function $f$.
        By the local version of Theorem~\ref{thm:alexander}, the family $\mathcal{F}$ is normal in some neighborhood of $0$ in $\mathbb{B}$.

        Consequently, some subsequence $\{f_{m_j}\}$ converges locally uniformly in a neighborhood of $0$ to a holomorphic function $f$.
        For every multi-index $\alpha$, the coefficient of $z^\alpha$ in $f_{m_j}$ is equal to $c_\alpha$ for all sufficiently large $j$.
        Since $f_{m_j}\to f$ locally uniformly, the Taylor coefficients of $f_{m_j}$ at $0$ converge to those of $f$.
        Hence the coefficient of $z^\alpha$ in the Taylor expansion of $f$ is $c_\alpha$.

        Therefore
        \[
            f(z)=\sum_\alpha c_\alpha z^\alpha
        \]
        in a neighborhood of $0$, and thus $F$ is convergent.
    \end{proof}

    \section{Main Result}

    The origins of the concept of the Kobayashi metric lie in Schwarz's lemma in complex analysis.

    \begin{defn}
        The infinitesimal Kobayashi metric on $\Omega$ is the function $F_K:\Omega\times\C^n\to\mathbb{R}^{+}$ defined by
        \[
            F_K(z,v) = \inf\set{ \alpha>0 \mid \exists\,\varphi\in\Omega(\D) \text{ with } \varphi(0)=z,\, \varphi'(0)=v/\alpha }\,.
        \]
    \end{defn}

    It is well known that, on the unit ball $\mathbb{B}$, the Bergman metric and the infinitesimal Kobayashi metric coincide up to a normalization constant; with the normalized Bergman metric, they coincide exactly (see, e.g.,~\cite{MR1098711,MR2446682}).

    \begin{defn}
        \label{defn-norm-function}
        Let $\Omega\subset\C^n$ be a domain and let $f:\Omega\to\C$ be a holomorphic function.
        We say that $f$ is normal if there exists a constant $C>0$ such that, for all $z\in\Omega$ and all $v\in\C^n$,
        \[
            L_z\left(\log(1+|f|^2),v\right) \leq C\,F_K(z,v)^2\,.
        \]
    \end{defn}

    \begin{thm}
        Let $\Omega$ be a domain in $\C^n$ such that the infinitesimal Kobayashi metric $F_K$ is nondegenerate.
        Let $f:\Omega\to\C$ be a holomorphic function.
        Then $f$ is normal on $\Omega$ if and only if the family
        \[
            \mathcal{F}= \set{f\circ\varphi \mid \varphi:\D\to\Omega\ \text{holomorphic}}
        \]
        is normal on $\D$.
    \end{thm}

    \begin{proof}
        Suppose $f$ is normal on $\Omega$.

        For any $\varphi\in\Omega(\D)$, put $g=f\circ\varphi$.
        Then
        \[
            \begin{aligned}
                g^\sharp(\lambda)^2 & = L_{\varphi(\lambda)}\left(\log(1+|f|^2),\varphi'(\lambda)\right) \\
                                        & \leq C\,F_K(\varphi(\lambda),\varphi'(\lambda))^2                \\
                                        & \leq C\left( \frac{2}{1-|\lambda|^2}\right)^2\,,
            \end{aligned}
        \]
        where the last inequality follows from the distance-decreasing property of the Kobayashi metric.
        Hence
        \[
            g^\sharp(\lambda) \leq \frac{2\sqrt{C}}{1-|\lambda|^2}, \qquad \lambda\in\D\,,
        \]
        uniformly for $g\in\mathcal{F}$.
        By Theorem~\ref{thm:criterion-norm}, $\mathcal{F}$ is a normal family on $\D$.

        Conversely, suppose that $\mathcal{F}$ is normal on $\D$.
        By Marty's theorem (Theorem~\ref{thm:marty:1}), there exists a constant $C>0$ such that
        \[
            L_{\varphi(0)}\left( \log(1+|f|^2),\varphi'(0) \right) = \frac{ |f'(\varphi(0))\varphi'(0)|^2 }{ (1+|f(\varphi(0))|^2)^2 }\leq C
        \]
        for every $\varphi\in\Omega(\D)$.

        Fix $z\in\Omega$ and $v\in\C^n$.
        If $v=0$, the desired inequality is trivial.
        Suppose that $v\neq0$.
        Since $F_K$ is nondegenerate, $F_K(z,v)>0$.
        Fix $\epsilon>0$, by the definition of $F_K$, there exist $\varphi\in\Omega(\D)$ such that
        \[
            \varphi(0)=z,\qquad \varphi'(0)=\frac{v}{\alpha}, \qquad \alpha<F_K(z,v)+\epsilon\,.
        \]

        Consequently,
        \[
            \begin{aligned}
                L_z\left(\log(1+|f|^2),v\right) & = \alpha^2 L_z\left( \log(1+|f|^2),\varphi'(0) \right) \\
                                                    & \leq C\alpha^2                                             \\
                                                    & \leq C(F_K(z,v)+\epsilon)^2\,.
            \end{aligned}
        \]

        Letting $\epsilon\to 0$ gives
        \[
            L_z\left(\log(1+|f|^2),v\right) \leq C\,F_K(z,v)^2\,.
        \]

        Thus $f$ is normal on $\Omega$.
    \end{proof}

    \section{Concluding Remark}

    Smoothly bounded domains in $\C^n$ generically have no holomorphic automorphisms; so our definition of normality is different from the one given above on $\mathbb{B}$.
    On the other hand, Definitions~\ref{defn-norm-function} and results in $\C^n$ are the natural generalizations of the result in $\C^1$.
    Some theorems related to those in the present work may be found in the book~\cite{MR4935535}.

    \section{Acknowledgments}

    We are grateful to the reviewers for their careful and constructive comments.
    We believe that their suggestions have significantly improved the correctness, clarity, and presentation of the manuscript.

    \bibliography{index}
\end{document}